\documentclass[12pt,a4paper]{article}
\usepackage{defns}
\usepackage[colour]{ajr}
\usepackage{url}
\usepackage{amsfonts}

\newcommand{\uj}{u_{j}}
\newcommand{\ibc}{\textsc{ibc}}

\newcommand{\difff}{{\cal D}}
\newcommand{\advvv}{{\cal C}}
\newcommand{\corrr}{{\cal B}}
\newcommand{\ddf}{\mathbb{D}}
\newcommand{\mudf}{\mathbb{C}}

\newcommand{\fdif}{\vec f_d}
\newcommand{\fadv}{\vec f_c}
\newcommand{\gadv}{\vec g_c}
\newcommand{\fcor}{\vec f_b}
\newcommand{\tr}[1]{{#1}^{\sf T}}

\title{Derive boundary conditions for holistic discretisations of 
Burgers' equation} 

\author{A. J. Roberts\thanks{Dept.\ Maths \& Comput, University of 
Southern Queensland, Toowoomba, Queensland 4352, \textsc{Australia}.  
\protect\url{mailto:aroberts@usq.edu.au}}}

\begin{document}

\maketitle

\begin{abstract}
I previously used Burgers' equation to introduce a new method of 
numerical discretisation of \pde{}s.
The analysis is based upon centre manifold theory so we are assured 
that the discretisation accurately models all the processes and their 
subgrid scale interactions.
Here I show how boundaries to the physical domain may be naturally 
incorporated into the numerical modelling of Burgers' equation.
We investigate Neumann and Dirichlet boundary conditions.
As well as modelling the nonlinear advection, the method naturally 
derives symmetric matrices with constant bandwidth to correspond to 
the self-adjoint diffusion operator.
The techniques developed here may be used to accurately model the 
nonlinear evolution of quite general spatio-temporal dynamical systems 
on bounded domains.
\end{abstract}

\tableofcontents

\section{Introduction}

We discuss the holistic discretisation of boundary conditions for 
partial differential equations.
Holistic discretisation \cite{Roberts98a} is based upon the support 
centre manifold theory gives to the nonlinear dynamics on finite grid 
spacing.
We expect such discretisation will have good stability and high 
accuracy on coarse grids because it systematically accounts for 
subgrid scale interactions.
To date we have considered periodic problems \cite{Roberts98a, 
MacKenzie00b, MacKenzie00a} and their initial conditions 
\cite{Roberts01a}.
Here we show how to incorporate different boundary conditions into the 
analysis.

As an illustrative example we restrict attention to the 
one-dimensional spatial discretisation of Burgers' equation
\begin{equation}
	\D tu+u\D xu=\DD xu\,,
	\label{eq:burg}
\end{equation}
which contains the important mechanisms of diffusion and nonlinear 
advection.
As example boundary conditions we consider the important cases of 
Dirichlet, \S\ref{Sdiri}, and Neumann boundary conditions, 
\S\ref{Sneum}.
The same techniques may be easily extended to other partial 
differential equations, other boundary conditions and higher spatial 
dimensions.

The centre manifold analysis is based upon dividing the domain into 
finite sized elements, each separated from their two neighbours by 
specially crafted artificial internal boundary conditions (\ibc{}s).
The form of the \ibc{}s~(\ref{eq:lbc}--\ref{eq:rbc}) generates a 
discretisation in the interior of the domain (\S\ref{Sadvec}) which is 
not only linearly consistent, as proved in~\cite{Roberts00a}, but 
which appears to be nonlinearly consistent to high order.
This observation is based upon this analysis of Burgers' equation and 
work in progress on the Kuramoto-Sivashinsky equation; further 
research is needed to prove it in general.

But the focus here is on the discretisation near the boundary of the 
domain.
Boundary conditions are easily incorporated simply by replacing an 
\ibc{} of an end element by a variant of the actual boundary 
condition: the Dirichlet boundary condition of fixed field~$u$ is 
implemented as~(\ref{eq:idiri}); the Neumann boundary condition of 
fixed flux~$u_x$ is~(\ref{eq:midn}).
The computer algebra of \S\ref{Aca} readily computes the effect these 
boundary conditions have on the discretisation near the boundary.
Novel features of this analysis are the following:
\begin{itemize}
	\item theoretical support is provided on finite grid spacing~$h$ as 
	explained in~\cite[e.g.]{Roberts98a,Roberts00a};
	
	\item the discretisation has a consistent bandwidth across the whole 
	domain including near the boundaries;

	\item high order discretisations are expressed in terms of the 
	matrices of the basic, second order, centred difference 
	operators, leading to appropriately symmetric discretisations of the 
	diffusion operator;

	\item time variations of boundary values, say~$a$, not only has a 
	direct effect but also involves time derivative factors, such 
	as~$\dot a$, which are more important on coarser grids due to the 
	time it takes changes in a boundary value to diffuse into the 
	element.

\end{itemize}

\section{Holistic discretisation of Burgers' equation in the 
interior}
\label{Sadvec}

We consider Burgers' equation~(\ref{eq:burg}) on some finite domain 
in~$x$.
Place grid points $x_j$ equi-spaced across the domain with constant 
spacing~$h$, and correspondingly define $\uj=u(x_j,t)$, that is, $\uj$ 
is the evolving field~$u$ evaluated at the $j$th grid point.
Divide the domain into $m$~elements with a grid point at the centre of 
each.
Following~\cite{Roberts00a} to apply centre manifold theory we 
distinguish each element using internal boundary conditions (\ibc{}s) 
in the discrete form
\begin{equation}
    \mu_x\delta_x v_j(x,t)=\gamma\mu\delta\,\uj
    \quad\mbox{and}\quad
    \delta_x^2 v_j(x,t)=\gamma\delta^2\uj\,,
    \label{eq:jbc}
\end{equation}
evaluated at $x=x_j$ where $v_j$ denotes the field in the $j$th 
element.
As shown in~\cite{Roberts00a}, this particular choice of \ibc{}s 
ensures that the resultant finite difference scheme is consistent to 
high order in~$h$ as the grid size~$h\to 0$\,.
In this section we repeat the construction of the holistic 
discretisation, following~\cite{Roberts98a}, away from the domain 
boundaries via centre manifold theory but with a new and convenient 
form of these discrete \ibc{}s.

In actually developing finite difference models the \ibc{}s may take 
any of many equivalent forms \cite[e.g.]{Roberts98a}.
In later sections we investigate the discretisation near a boundary of 
the domain: the element adjacent a boundary has one real boundary 
and one artificial internal boundary.
Thus it is appropriate to rewrite the two \ibc{}s in~(\ref{eq:jbc}) in 
the form of two conditions, one at the left edge of each internal 
element, and one at the right.
Recall that the difference operators $\mu\pm\half\delta=E^{\pm 1/2}$ 
\cite[p65, e.g.]{npl61} so that to the first \ibc\ in~(\ref{eq:jbc}) 
add\slash subtract half the second to give the equivalent \ibc{}s
\begin{eqnarray}
	\delta_x v_j(x,t)&=&\gamma\delta u_{j-1/2}
	\quad\mbox{at $x=x_{j-1/2}$}\,,
	\label{eq:lbc} \\ \mbox{and}\quad
	\delta_x v_j(x,t)&=&\gamma\delta u_{j+1/2}
	\quad\mbox{at $x=x_{j+1/2}$}\,.
	\label{eq:rbc} 
\end{eqnarray}
The \ibc~(\ref{eq:rbc}) is to hold at the right-hand side of each 
element and the \ibc~(\ref{eq:lbc}) is to hold at the left.
The introduced parameter~$\gamma$, when non-zero, couples each element 
to its neighbour.
When $\gamma=0$ these \ibc{}s effectively insulate each element from 
its neighbours and forms the basis of the centre manifold analysis; 
whereas when $\gamma=1$ they assert that the field in the $j$th 
element when extended into the surrounding elements has the same 
differences centred across each internal boundary as that given by the 
grid point values.
Thus evaluating the model at~$\gamma=1$ forms the relevant 
discretisation.
These \ibc{}s apply to all elements except for the leftmost and 
rightmost elements, the ones adjacent to the boundary: in the leftmost 
element the left-hand \ibc, (\ref{eq:lbc}) with $j=1$, is replaced by 
the actual boundary condition; in the rightmost element the right-hand 
\ibc, (\ref{eq:rbc}) with $j=m$, is replaced by the actual boundary 
condition.
In this paper I analyse the left boundary of the 
domain---discretisations near the right boundary are similar by 
symmetry.

In the interior of the domain the boundaries of the domain have no 
influence upon the discretisation.
Each element in the interior has \ibc{}s~(\ref{eq:lbc}) 
and~(\ref{eq:rbc}).
Executing the computer algebra program in Appendix~\ref{Aca}, adapted 
from~\cite{Roberts98a}, the subgrid structure of the solution 
field~$u(x,t)$ is, in terms of $\xi=(x-x_j)/h$,
\begin{eqnarray}
	v_j & = & \uj +\gamma\left[ \xi\mu\delta\uj +\half\xi^2\delta^2\uj 
	\right]
	+\gamma^2\left[ \sfrac16(\xi^3-\xi)\mu\delta^3\uj  
	 +\sfrac1{24}(\xi^4-\xi^2)\delta^4\uj \right]
	 \nonumber\\&&{}
	 +\gamma h\sfrac16(\xi^3-\xi)\uj\delta^2\uj
	 +\Ord{\|\vec u\|^3+\gamma^3}\,.
	\label{eq:cm}
\end{eqnarray}
The corresponding evolution on this centre manifold forms the holistic 
discretisation in the interior:
\begin{eqnarray}
	\dot\uj & = & \phantom{{}+}\frac{1}{h^2}\left[ \gamma\delta^2\uj 
	-\frac{\gamma^2}{12}\delta^4\uj +\frac{\gamma^3}{90}\delta^6\uj 
	\right]
	\nonumber  \\
	 &  & {}-\frac{1}{h}\uj\left[ \gamma\mu\delta\uj 
	 -\frac{\gamma^2}{6}\mu\delta^3\uj \right]
	 +\frac{\gamma^2}{24h}( \delta^2\uj\,\mu\delta^3\uj 
	 +\delta^4\uj\,\mu\delta\uj )
	\nonumber  \\
	 &  & {}+\frac{\gamma}{12}\uj^2\delta^2\uj
	 +\Ord{\|\vec u\|^4+\gamma^4}\,.
	\label{eq:hmod}
\end{eqnarray}
As discussed previously~\cite{Roberts98a} when the coupling 
parameter~$\gamma=1$: the first line gives successive approximations 
to the diffusion term~$u_{xx}$; the second line gives approximations 
to the nonlinear advection; and lastly, the cubic nonlinear term, 
$\uj^2\delta^2\uj$ in the last line above, accounts for subgrid scale 
interactions between the advection and diffusion and acts to stabilise 
the numerical model.
Higher-order terms are easily found by the computer algebra of 
Appendix~\ref{Aca} but for clarity are not presented here.
In the limit as the grid spacing~$h\to 0$ higher-order discretisations 
have the equivalent \pde{}
\begin{eqnarray}
	\D tu & = & \gamma\left[ -u\D xu+\DD xu \right]
	+\frac{h^2}{12}\gamma(1-\gamma)\left[ u_{xxxx} -2uu_{xxx} 
	 +u^2u_{xx} \right]
	\nonumber  \\
	 &  & {}+\frac{h^4}{720}\gamma(1-\gamma)
	 \left[ 2\gamma(-5u_x^2u_{xx} -9uu_{xx}^2 -25uu_xu_{xxx} 
	 +15u_{xx}u_{xxx} 
	 \right.\nonumber\\&&\left.\quad{}
	 +15u_xu_{xxxx} -2u^2u_{xxxx} ) 
	 +(1-4\gamma)( 2\partial_x^6u -6u\partial_x^5u +5u^2\partial_x^4u ) 
	 \right]
	\nonumber  \\
	 &  & {}+\Ord{\|u\|^4,\gamma^4,h^6}\,.
	\label{eq:equivde}
\end{eqnarray}
See that upon substituting $\gamma=1$ to recover a discretisation for 
the Burgers' equation~(\ref{eq:burg}), we would find an equivalent 
\pde{} to an error $\Ord{h^6,\|u\|^4}$.
Analogous lower order discretisations are obtained, as promised by the 
analysis in~\cite{Roberts00a}, when we truncate the 
discretisation~(\ref{eq:hmod}) to lower orders in the coupling 
parameter~$\gamma$.
But observe the new feature that the discrete 
form~(\ref{eq:lbc}--\ref{eq:rbc}) of the \ibc{}s lead to 
discretisations which are not only linearly consistent, as assured 
in~\cite{Roberts00a}, but are also nonlinearly consistent.
Further research is needed to establish nonlinear consistency in 
general.

\section{Dirichlet boundary conditions applied at a grid point}
\label{Sdiri}

Consider the case of Dirichlet conditions on the boundary of the 
domain of prescribed~$u$ at a grid point; without loss of generality 
say
\begin{equation}
	u=a(t)\quad\mbox{at $x=x_0$}\,.
	\label{eq:diri}
\end{equation}
This boundary condition is included in the analysis simply by 
replacing the left-hand \ibc\ in the leftmost element, (\ref{eq:lbc}) 
with $j=1$, by (as if $u_0=a$ in~(\ref{eq:lbc}))
\begin{equation}
		\delta_x v_1(x,t)=\gamma(u_1-a)
	\quad\mbox{at $x=x_{1/2}$}\,.
	\label{eq:idiri}
\end{equation}
(Implicitly the first element then extends from $x_0=x_1-h$ 
to~$x_1+h/2$.) When the coupling parameter~$\gamma=0$ this \ibc\ 
effectively insulates the first element from the conditions at the 
domain boundary.
However, when $\gamma=1$, since $v_j(x_j,t)=\uj$, this reduces 
to~(\ref{eq:diri}) by requiring $v_1(x_0,t)=a$\,.
Hence, the centre manifold derivation is based on $\gamma=0$ as 
explained in~\cite{Roberts98a} and evaluated at $\gamma=1$ to obtain a 
discretisation of Burgers' equation.

We then solve for the subgrid fields in the elements near the 
boundary, $j=1,2,\ldots$\,.
The computer algebra program in Appendix~\ref{Aca} implements this 
boundary condition when \verb|dirichlet:=1|\,.
The influence of this specified boundary value affects a number of 
elements near the boundary equal to the order of~$\gamma$ retained in 
the analysis, nothing else: in the interior the discretisation 
is~(\ref{eq:hmod}); whereas for elements the near the boundary and for 
errors $\Ord{\|\vec u\|^4+\gamma^4}$ we find the evolution to be of 
the form
\begin{eqnarray}
	\left[\begin{array}{c}
		\dot u_1  \\
		\dot u_2  \\
		\dot u_3
	\end{array}\right] &=& 
	\frac{1}{h^2}\left(\difff \vec u+\fdif\right)
	-\frac{1}{h}\left(U\advvv \vec u+\gadv(\vec u)+\fadv\right)
	\nonumber\\&&{}
	+\left(U^2\corrr \vec u+\fcor\right)
	+\Ord{\|\vec u\|^4+\gamma^4}\,,
	\label{eq:dirix}
\end{eqnarray}
where $U=\mbox{diag}(u_1,u_2,u_3)$ is the diagonal matrix of grid 
velocities and the three parts of the right-hand side represent 
respectively the discretisation of the diffusion, the nonlinear 
advection, and the leading order interaction between advection and 
diffusion.
These parts include the forcing due to the time dependent boundary 
value~$a$.
Here the various terms are found to be:
\begin{eqnarray}
	\difff&=&
	\gamma\left[\begin{array}{rrrrr}
		-2 & 1 &  &   \\
		1 & -2 & 1 &  &\cdots \\
		 & 1 & -2 & 1
	\end{array}\right]
	-\frac{\gamma^2}{12}\left[\begin{array}{rrrrrr}
	 	5 & -4 & 1 &  &   \\
	 	-4 & 6 & -4 & 1 &  &\cdots \\
	 	1 & -4 & 6 & -4 & 1
	 \end{array}\right]
	 \nonumber\\&&{}
	 +\frac{\gamma^3}{90}\left[\begin{array}{rrrrrrr}
	 	-14 & 14 & -6 & 1 &  &   \\
	 	14 & -20 & 15 & -6 & 1 &  &\cdots \\
	 	-6 & 15 & -20 & 15 & -6 & 1
	 \end{array}\right]
	 +\Ord{\gamma^4}
	\label{eq:difff}\\
	\advvv&=&
	\frac{\gamma}{2}\left[\begin{array}{rrrrr}
		0 & 1 &  &   \\
		-1 & 0 & 1 &  &\cdots \\
		 & -1 & 0 & 1
	\end{array}\right]
	 -\frac{\gamma^2}{12}\left[\begin{array}{rrrrrr}
	 	0 & -2 & 1 &  &   \\
	 	2 & 0 & -2 & 1 &  &\cdots \\
	 	-1 & 2 & 0 & -2 & 1
	 \end{array}\right]
	 \nonumber\\&&{}
	 +\Ord{\gamma^3}
	\label{eq:advvv}\\
	\gadv&=&
	\frac{\gamma^2}{24}\left[\begin{array}{c}
		\half \tr{\vec u}G_1\vec u  \\
		\half \tr{\vec u}G_2\vec u  \\
		\half \tr{\vec u}G_3\vec u
	\end{array}\right]
	\quad\mbox{where}\quad
	G_1=\left[\begin{array}{ccc}
		2 & 5 & -1  \\
		5 & -6 & 1  \\
		-1 & 1 & 0
	\end{array}\right]
	\nonumber\\&&\!\!\!
	G_2=\left[\begin{array}{cccc}
		6 & -5 & 0 & 0  \\
		-5 & 0 & 5 & -1  \\
		0 & 5 & -6 & 1  \\
		0 & -1 & 1 & 0
	\end{array}\right]\ 
	G_3=\left[\begin{array}{ccccc}
		0 & -1 & 1 & 0 & 0  \\
		-1 & 6 & -5 & 0 & 0  \\
		1 & -5 & 0 & 5 & -1  \\
		0 & 0 & 5 & -6 & 1  \\
		0 & 0 & -1 & 1 & 0
	\end{array}\right]
	\label{eq:gs}
\end{eqnarray}
and $\corrr=\frac{1}{12}\difff+\Ord{\gamma}$.
Denote by $\ddf$ the matrix appearing above in $\difff$ linear 
in~$\gamma$ and then row extended across the interior of the domain: 
$\ddf$ is the matrix of the second-order centred approximation to the 
second derivative, $\delta^2$.
Observe that the order~$\gamma^2$ and~$\gamma^3$ matrices in $\difff$ 
simply correspond to~$\ddf^2$ and~$\ddf^3$.
Thus the discretisation of the diffusion term, across the entire 
domain including the near boundary elements, is simply $\gamma\delta^2 
-\frac{\gamma^2}{12}\delta^4 +\frac{\gamma^3}{90}\delta^6$ seen in the 
first line of the interior discretisation,~(\ref{eq:hmod}), but with 
the matrix~$\ddf$ replacing the centred difference~$\delta^2$.
Thus this approach generates an appropriately symmetric discretisation 
of the self-adjoint diffusion term~$u_{xx}$.
The nonlinear stabilisation term, $U^2\corrr\vec u$, also corresponds 
to replacing $\delta^2$ in~(\ref{eq:hmod}) by~$\ddf$.
Similarly, denote by $\mudf$ the matrix appearing above in $\advvv$ 
linear in~$\gamma$ (including the factor $\half$) and then row 
extended across the interior of the domain: $\mudf$ is the matrix of 
the second-order centred approximation to the first derivative, 
$\mu\delta$.
Observe that the order~$\gamma^2$ matrix in $\advvv$ is~$(\ddf\mudf 
+\mudf\ddf)/2$ corresponding to the average of different permutations 
of the centred difference operators $\mu\delta^3$.
Thus again the discretisation of the advection terms across the entire 
domain is the interior discretisation, $\gamma\mu\delta 
-\frac{\gamma^2}{6} (\mu\delta\delta^2 +\delta^2\mu\delta)/2$, with 
the truncated~$\mudf$ and~$\ddf$ replacing the difference 
operators~$\mu\delta$ and~$\delta^2$ respectively.
Additional terms represented by $\gadv$ for the discretisation of the 
advection correspond to the nonlinear higher order term 
$\frac{\gamma^2}{24h}( \delta^2\uj\,\mu\delta^3\uj 
+\delta^4\uj\,\mu\delta\uj )$ in the model~(\ref{eq:hmod}).
These are pleasing patterns.

Also see that by truncating to a fixed power of the inter-element 
coupling parameter~$\gamma$ we will obtain a discretisation that has 
constant bandwidth across the whole domain.
This constant bandwidth will always be derived in this holistic 
approach, see for another example the Neumann boundary conditions in 
the next section, because the truncation at a fixed power of the 
inter-element coupling~$\gamma$ controls how many neighbouring 
elements interact with any given element.
Although there is a lower order (in~$h$) of consistency near the 
domain boundaries, the support by centre manifold theory is the same 
across the whole domain.
As discussed in~\cite{Roberts98a}, the reason for this support is that 
the theory applies to the solutions of Burgers' \pde~(\ref{eq:burg}) 
in the entire domain, not just in some locale.

Now consider the forcing from the boundary.
Using $\vec a=(a,h^2\dot a)$, it is
\begin{eqnarray}
	\fdif&=&
	\left[\begin{array}{rr}
		\gamma+\frac{\gamma^2}{6}+\frac{\gamma^3}{18}
		& {}-\frac{\gamma}{12}-\frac{\gamma^2}{45}-\frac{\gamma^3}{112}
		\\[1ex]
		-\frac{\gamma^2}{12}-\frac{2\gamma^3}{45}
		& {}+\frac{\gamma^2}{90}+\frac{\gamma^3}{140}  \\[1ex]
		\frac{\gamma^3}{90}
		& {}-\frac{\gamma^3}{560}
	\end{array}\right]\vec a
	 +\Ord{\gamma^4,\ddot a}\,,
	\label{eq:fdif}\\
	\fadv&=&
	\gamma\left[\begin{array}{cc}
		\half u_1 & -\frac{1}{24}u_1\\
		0&0\\
		0&0
	\end{array}\right]\vec a
	\nonumber\\&&{}
	+\frac{\gamma^2}{24}\left[\begin{array}{cc}
		4u_1+3a & \frac{1}{15}(7u_1-2u_2-9a)+\frac{5}{168}h^2\dot a\\
		-(u_1+u_2)& \frac{2}{15}(u_1+u_2)\\
		0&0
	\end{array}\right]\vec a
	\nonumber\\&&{}
	 +\Ord{\gamma^3,\ddot a}\,,
	\label{eq:fadv}\\
	\fcor&=&
	\frac{\gamma}{12}\left[\begin{array}{cc}
		u_1^2
		& -\frac{1}{15}u_1^2 \\
		0&0\\
		0&0
	\end{array}\right]\vec a
	 +\Ord{\gamma^2,\ddot a}\,.
	\label{eq:fcor}
\end{eqnarray}
See that this holistic approach not only involves the value~$a$ of the 
field at the boundary in the diffusion term, it also involves~$a$ in 
the nonlinear advection (in~$\fadv$) and in the nonlinear 
stabilisation (in~$\fcor$).
But further, it also involves time derivatives of the boundary 
value~$a$.
The reason is clearly that changes in~$a$ take time to advect and 
diffuse into the interior of the first few elements and so the effect 
of changes in the boundary value~$a$ upon the evolution of the grid 
values will lag, as seen by the opposite sign of the coefficients in 
the two columns of~(\ref{eq:fdif}).
This lag increases with the element size~$h$ and hence accounts for 
the $h^2$~factor multiplying every~$\dot a$.
Similarly, though not recorded above, higher order analysis shows each 
second derivative appears only in the combination~$h^4\ddot a$.
Such effects are important when we try to use, in large scale 
problems, the expected accuracy and stability of these holistic 
discretisations on coarse grids.

\section{Apply Neumann boundary conditions at a midpoint}
\label{Sneum}

Consider \pde{}s with Neumann boundary conditions of prescribed 
gradient of the field~$u$.
Two different sorts of numerical approximations are commonly developed 
for such a boundary condition: a grid point is placed at the boundary 
(as for the Dirichlet problem of \S\ref{Sdiri}); or the boundary is 
arranged to be midway between two grid points.
For the first case I found that the discretisation of the diffusion is 
expressed in terms of non-symmetric matrices.
For the second case, the matrices are symmetric which again 
corresponds to the self-adjoint nature of the diffusion operator.
Thus I report here on the second case where the boundary is at a 
midpoint of the grid.

Without loss of generality, let the Neumann boundary condition of 
prescribed gradient at a grid midpoint be\footnote{Neumann conditions 
\emph{at a grid point}, say $x_1$, may be incorporated into the 
analysis by similarly requiring $h\partial u/\partial x=\gamma a$ 
at~$x=x_1$ in place of the left-hand \ibc~(\ref{eq:lbc}).}
\begin{equation}
	h\D xu=\gamma{a(t)}\quad\mbox{at $x=x_{1/2}$}\,,
	\label{eq:midn}
\end{equation}
where, as for the \ibc{}s, we actually are interested in the case 
$\gamma=1$.
To supplement this left-hand boundary condition on the leftmost 
element we use the \ibc~(\ref{eq:rbc}) as before.
We execute the computer algebra program of Appendix~\ref{Aca} with 
options \verb|dirichlet:=0| and \verb|midpoint:=1|\,.\footnote{The 
computer algebra algorithm takes a few more iterations to complete 
because the algorithm is tuned to the discrete form of the 
\ibc~(\ref{eq:lbc}) whose left-hand side is only approximately that 
of~(\ref{eq:midn}).} Again the interior discretisation 
is~(\ref{eq:hmod}) whereas in elements the near the boundary we find 
the grid values evolve according to the form~(\ref{eq:dirix}) where 
now the matrices are
\begin{eqnarray}
	\difff&=&
	\gamma\left[\begin{array}{rrrrr}
		-1 & 1 &  &   \\
		1 & -2 & 1 &  &\cdots \\
		 & 1 & -2 & 1
	\end{array}\right]
	-\frac{\gamma^2}{12}\left[\begin{array}{rrrrrr}
	 	2 & -3 & 1 &  &   \\
	 	-3 & 6 & -4 & 1 &  &\cdots \\
	 	1 & -4 & 6 & -4 & 1
	 \end{array}\right]
	 \nonumber\\&&{}
	 +\frac{\gamma^3}{90}\left[\begin{array}{rrrrrrr}
	 	-5 & 9 & -5 & 1 &  &   \\
	 	9 & -19 & 15 & -6 & 1 &  &\cdots \\
	 	-5 & 15 & -20 & 15 & -6 & 1
	 \end{array}\right]
	 +\Ord{\gamma^4}\,,
	\label{eq:mdifff}\\
	\advvv&=&
	\frac{\gamma}{2}\left[\begin{array}{rrrrr}
		-1 & 1 &  &   \\
		-1 & 0 & 1 &  &\cdots \\
		 & -1 & 0 & 1
	\end{array}\right]
	 -\frac{\gamma^2}{12}\left[\begin{array}{rrrrrr}
	 	1 & -2 & 1 &  &   \\
	 	1 & 0 & -2 & 1 &  &\cdots \\
	 	-1 & 2 & 0 & -2 & 1
	 \end{array}\right]
	 \nonumber\\&&{}
	 +\Ord{\gamma^3}\,,
	 \\
	\gadv&=& \frac{\gamma}{24}\left[\begin{array}{c}
		-u_1^2+u_2^2  \\
		0  \\
		0
	\end{array}\right]
	+\frac{\gamma^2}{24}\left[\begin{array}{c}
		\half \tr{\vec u}G_1\vec u  \\
		\half \tr{\vec u}G_2\vec u  \\
		\half \tr{\vec u}G_3\vec u
	\end{array}\right]
	\quad\mbox{where}\nonumber\\&&
	G_1=\left[\begin{array}{ccc}
		-\frac{49}{20} & \frac{19}{15} & -\frac{11}{15}  \\
		\frac{19}{15} & -\frac{341}{60} & 1  \\
		-\frac{11}{15} & 1 & 0
	\end{array}\right]
	\quad
	G_2=\left[\begin{array}{cccc}
		\frac{23}{6} & -\frac{47}{12} & 0 & 0  \\
		-\frac{47}{12} & 0 & 5 & -1  \\
		0 & 5 & -6 & 1  \\
		0 & -1 & 1 & 0
	\end{array}\right]
	\label{eq:madvvv}\\
	\corrr&=&\frac{\gamma}{12}\tilde\ddf+
	\frac{\gamma}{12}\left[\begin{array}{rrrrr}
		-\frac{1}{48} & \frac{1}{48} &  &   \\
		0 & 0 & 0 &  &\cdots \\
		 & 0 & 0 & 0
	\end{array}\right]
	 +\Ord{\gamma^2}\,,
	\label{eq:mcorrr}
\end{eqnarray}
and $G_3$ is as in~(\ref{eq:gs}).
Denote by $\tilde\ddf$ the symmetric matrix appearing above in 
$\difff$ linear in~$\gamma$ and then row extended across the interior 
of the domain: $\tilde\ddf$ is the matrix of the second-order centred 
approximation to the second derivative, $\delta^2$, but incorporating 
a zero derivative condition to the left of the grid.
Observe in the order~$\gamma^2$ and~$\gamma^3$ matrices of~$\difff$ 
that the discretisation near the boundary of the diffusion term, 
(\ref{eq:mdifff}), is simply the interior discretisation seen in the 
first line of~(\ref{eq:hmod}) with the truncated~$\tilde\ddf$ 
replacing the centred difference~$\delta^2$.
Although here there is no clear pattern in $\advvv$ nor~$\corrr$, as 
written above see that: $\corrr$ is numerically close to 
$\tilde\ddf/12$; and, upon denoting $\tilde\mudf$ as the matrix linear 
in $\gamma$ in~$\advvv$, the $\gamma^2$ matrix is 
$(\tilde\ddf\tilde\mudf+\tilde\mudf\tilde\ddf)/2$ just as for the 
Dirichlet boundary conditions.
However, in this case the identification is a little forced because 
numerically small discrepancies have been gathered into $\gadv$ along 
with indistinguishable similar terms corresponding to 
$\frac{\gamma^2}{24h}( \delta^2\uj\,\mu\delta^3\uj 
+\delta^4\uj\,\mu\delta\uj )$ appearing in the interior 
discretisation~(\ref{eq:hmod}).
Nonetheless the closeness of the match and the appropriate symmetry is 
reassuring.

The forcing from the boundary takes the form
\begin{eqnarray}
	&&\fdif=
	\left[\begin{array}{rr}
		-\gamma-\frac{\gamma^2}{12}-\frac{\gamma^3}{45}
		& {}+\frac{\gamma}{24}+\frac{11\gamma^2}{1440}
		+\frac{\gamma^3}{378}  \\[1ex]
		+\frac{\gamma^2}{12}+\frac{\gamma^3}{30} 
		& {}-\frac{11\gamma^2}{1440}-\frac{\gamma^3}{252}  \\[1ex]
		-\frac{\gamma^3}{90}
		& {}+\frac{\gamma^3}{756}
	\end{array}\right]\vec a
	 +\Ord{\gamma^4,\ddot a}\,,
	\label{eq:mfdif}\\
	&&\fadv=
	\gamma\left[\begin{array}{cc}
		-\frac{11}{24} u_1 & +\frac{31}{960}u_1\\
		0&0\\
		0&0
	\end{array}\right]\vec a
	\nonumber\\&&{}
	+\frac{\gamma^2}{24}\left[\begin{array}{cc}
		-\frac{1}{6}u_1-\frac{11}{60}u_2-\frac{199}{120}a 
		& \frac{924u_1-577u_2-3223a}{10080}-\frac{757}{48384}h^2\dot a\\
		-(\frac{11}{12}u_1+u_2)
		& \frac{53}{480}u_1+\frac{11}{120}u_2\\
		0&0
	\end{array}\right]\vec a
	\nonumber\\&&{}
	 +\Ord{\gamma^3,\ddot a}\,,
	\label{eq:mfadv}\\
	&&\fcor=
	\frac{\gamma}{12}\left[\begin{array}{cc}
		-\frac{49}{48}u_1^2
		& \frac{637}{5760}u_1^2 \\
		0&0\\
		0&0
	\end{array}\right]\vec a
	 +\Ord{\gamma^2,\ddot a}\,.
	\label{eq:mfcor}
\end{eqnarray}
As discussed in the previous section, the boundary value for the 
flux~$a$ appears in a wide range of terms in the discretisation near 
the boundary.
The reason again is that the flux feeds into the subgrid scale fields 
of the boundary elements and so affects the interaction between the 
various physical processes.
The scope for such interaction increases with increasing element 
size~$h$ and so accounts for the appearance of the $h^2$ factor in 
front of the time derivative~$\dot a$.
We need to know such effects on coarse grids.

\section{Conclusion}

We have considered the discretisation of Burgers' 
equation~(\ref{eq:burg}) on a finite domain as a worked example of 
incorporating physical boundary conditions into the derivation of 
holistic discretisations.  The two most common physical boundary 
conditions were considered in~\S\ref{Sdiri} and~\S\ref{Sneum}.  
Crucially we found: it is easy to maintain the symmetry appropriate to 
self-adjoint differential operators; discretisations are developed 
with constant bandwidth across the whole domain; and lastly our 
resolution of subgrid scale processes shows how time derivatives of 
the boundary forcing should also be included in the discretisation.

Additionally, using the discrete form of the inter-element boundary 
conditions~(\ref{eq:lbc}--\ref{eq:rbc}) we observe for the first time 
a high order consistency of the nonlinear dynamics of Burgers' 
equation.

\paragraph{Acknowledgements:} this research was supported in part by a 
grant from the Australian Research Council.

\bibliographystyle{plain}
\bibliography{ajr,bib,new}

\begin{thebibliography}{1}

\bibitem{MacKenzie00b}
T.~MacKenzie and A.~J. Roberts.
\newblock The dynamics of reaction diffusion equations lead to a holistic
  discretisation.
\newblock In R.~L. May, G.~F. Fitz-Gerald, and I.~H. Grundy, editors, {\em EMAC
  2000 Proceedings. Proceedings of the fourth biennial Engineering Mathematics
  and Applications Conference}, pages 199--202, 2000.

\bibitem{MacKenzie00a}
T.~Mackenzie and A.~J. Roberts.
\newblock Holistic finite differences accurately model the dynamics of the
  {Kuramoto-Sivashinsky} equation.
\newblock {\em ANZIAM~J.}, 42(E):C918--C935, 2000.
\newblock [Online] \url{http://anziamj.austms.org.au/V42/CTAC99/Mack}.

\bibitem{npl61}
{National Physical Laboratory}.
\newblock {\em Modern Computing Methods}, volume~16 of {\em Notes on Applied
  Science}.
\newblock Her Majesty's Stationary Office, 1961.

\bibitem{Roberts96a}
A.~J. Roberts.
\newblock Low-dimensional modelling of dynamics via computer algebra.
\newblock {\em Comput. Phys. Comm.}, 100:215--230, 1997.

\bibitem{Roberts00a}
A.~J. Roberts.
\newblock A holistic finite difference approach models linear dynamics
  consistently.
\newblock Technical report, [\url{http://arXiv.org/abs/math.NA/0003135}], March
  2000.

\bibitem{Roberts98a}
A.~J. Roberts.
\newblock Holistic discretisation ensures fidelity to {Burgers'} equation.
\newblock {\em Applied Numerical Modelling}, 37:371--396, 2001.

\bibitem{Roberts01a}
A.~J. Roberts.
\newblock Holistic projection of initial conditions onto a finite difference
  approximation.
\newblock Technical report, [\url{http://arXiv.org/abs/math.NA/0101205}], 2001.

\end{thebibliography}


\appendix
\section{Computer algebra derives different boundary discretisations}
\label{Aca}

Straightforward computer algebra programs are written to find the 
centre manifold and the evolution thereon \cite[e.g.]{Roberts96a}.
To ensure correctness and to provide a basis for further work I 
include the computer algebra code.
This replaces extensive recording of elementary algebraic steps in the 
derivation of the results.

I implement the construction algorithm in \textsc{reduce}\footnote{At 
the time of writing, information about \texttt{reduce} was available 
from Anthony C.~Hearn, RAND, Santa Monica, CA~90407-2138, USA. 
\url{mailto:reduce@rand.org}} The overall plan of the algorithm is to 
iteratively satisfy Burgers' equation~(\ref{eq:burg}) and the relevant 
internal~(\ref{eq:lbc}--\ref{eq:rbc}) and actual boundary conditions.  
A general interior element is analysed in lines~47--53 while the 
\verb|o:=3| elements near the boundary are analysed in the loop of 
lines~54--71.  Although there are many details in the program, the 
correctness of the results are \emph{only determined} by driving to 
zero (lines~53, 70 and~73) the residuals of Burgers' equation in each 
element and the internal and domain boundary conditions: lines~47--9 
evaluate the residuals for an arbitrary interior element; lines~55--8 
for near domain boundary elements; and lines~60 or~62 the domain 
boundary condition.  The other details, such as the updates computed 
in lines~50--2 and~67--9, only affect the rate of convergence to the 
ultimate answer.

{\small
\verbatimlisting{burg.red}
}

\end{document}